\documentclass[12pt]{article}
\usepackage{kuvio}

 \usepackage{amssymb}
 \usepackage{latexsym}

 \newcommand{\im}{\mathop{\rm im}\nolimits}
 \newcommand{\codim}{\mathop{\rm codim}\nolimits}
 \newcommand{\sk}{\mathop{\rm sk}\nolimits}
 \newcommand{\id}{\mathop{\rm Id}\nolimits}
 \newcommand{\R}{{\mathbb R}}
 \newcommand{\Z}{{\mathbb Z}}
 
 \newcommand{\eps}{{\varepsilon}}
 \newcommand{\definition}{\noindent{\bf Definition\ }}
 \newcommand{\example}{\noindent{\bf Example\ }}

 \newtheorem{theorem}{Theorem}
 \newtheorem{lemma}{Lemma}[section]
 \newtheorem{remark}{Remark}
 
 \newtheorem{proposition}{Proposition}

 \newcommand{\beginproof}{\begin{trivlist}
 \rm\item[\hspace{\labelsep}{ P\,r\,o\,o\,f. }]}
 \newcommand{\proofend}{$\Box$\end{trivlist}}

\author{Fedor S. Duzhin}
\title
{Lower bounds for the number of closed billiard trajectories
of period $2$ and $3$ in manifolds embedded in Euclidean space}
\date{November 2001}

\begin{document}
\maketitle

\section{Introduction}
Let $M$ be a smooth closed $m$-dimensional manifold
embedded in Euclidean space $\R^n$. 
An ordered set $(x_1,\dots,x_p)\in M\times\dots\times M$
is called {\it a periodic (closed) billiard trajectory} if 
the condition
$$
\frac{x_i-x_{i+1}}{\|x_i-x_{i+1}\|}+
\frac{x_i-x_{i-1}}{\|x_i-x_{i-1}\|}\perp T_{x_i}M
$$
holds for each $i\in\Z_p$.

In this paper we find lower bounds for the number of
closed billiard trajectories of
periods $2$ and $3$ in an arbitrary manifold.

Let us introduce some notations:
$$
M^{\times p}=M\times\dots\times M
$$
is $p$th Cartesian power of the space $M$,
$$
\tilde{\Delta}^{(k)}=\bigcup_{i\in\Z_p}
\left\{\left(x_1,\dots,x_p\right):
x_i=x_{i+1}=\dots=x_{i+p-k-1}\right\}\subset M^{\times p}
$$
is the $k$th diagonal of $M^{\times p}$. Clearly, we have
$$
\left\{(x,x,\dots,x)\right\}=
\tilde{\Delta}^{(0)}\subset\tilde{\Delta}^{(1)}\subset
\dots\subset\tilde{\Delta}^{(p-2)}\subset\tilde{\Delta}^{(p-1)}
=M^{\times p}.
$$
By $\tilde\Delta$ we denote 
$\tilde{\Delta}^{(p-2)}=\bigcup_{i\in\Z_p}\{x_i=x_{i+1}\}
\subset M^{\times p}$. Let
$$
\tilde{l}=\sum_{i\in\Z_p}\|x_i-x_{i+1}\|:M^{\times p}\to\R
$$
be the length function of a closed polygon. The
function $\tilde{l}$ is smooth outside of the diagonal 
$\tilde{\Delta}$. The symmetric group $S_p$ acts on
$M^{\times p}$ by permutations of factors, and so does
the dihedral group $D_p\subset S_p$.
Let
$$
D^p(M)=M^{\times p}/D_p,\ 
\Delta^{(k)}=\tilde{\Delta}^{(k)}/D_p\subset D^p(M),\ 
l=\tilde{l}/D_p:D^p(M)\to\R.
$$
By $RD^p(M)$ we denote the reduced dihedric power of $M$,
that is $D^p(M)$ with the diagonal 
$\Delta=\Delta^{(p-2)}$ contracted to a point.

It is well known that $p$-periodical billiard trajectories
are the critical points of the function $f$ outside of
$\Delta$.
Thus our problem is to estimate the minimal number $BT_p(M)$ of
critical points of the function $l$ for all generic embeddings
of the manifold $M$ in Euclidean space.
Here we say that an embedding $M\to\R^n$ is generic
if the function $l$ has only non-degenerate critical points
outside of $\Delta$.

Recall the following facts:
\medskip
\begin{theorem} (see \cite{duzhin})
Let $p$ be a prime integer. Then
$$
BT_p(M)\ge\sum\dim H_q(M^{\times p}/D_p,\Delta;\Z_2).
$$
\end{theorem}
\medskip

\begin{theorem} (M.~Farber, S.~Tabachnikov, \cite{farber}, \cite{farber2})
Let $p$ be an odd prime. Then
$$
\sum\dim\tilde{H}_q(RD^p(S^m);\Z_2)=m(p-1).
$$
\end{theorem}

\bigskip 
In this paper we estimate
the number of closed billiard trajectories
of period $2$ and $3$ for an arbitrary manifold. 
The estimate
for period $2$ was proved by P.~Pushkar
in \cite{pushkar}, we give a new proof.

\begin{theorem}
Let $M$ be a smooth closed $m$-dimensional
manifold. By $B$ denote the sum of Betti numbers
$\sum_{q=0}^m\dim H_q(M;\Z_2)$. Then we
have the following estimates for the number
of periodical billiard trajectories:
$$
BT_2(M)\ge\frac{B^2+(m-1)B}{2},
$$
$$
BT_3(M)\ge\frac{B^3+3(m-1)B^2+2B}{6}.
$$
\end{theorem}

\section{A.~Dold's theory of homology of symmetric products}
In this section we consider $\Z_2$-modules.
\medskip

\definition
{\it An FD-module} $K$ is a sequence of modules
$K_q$, $q=0,1,\dots$ and module morphisms
$\partial^i_q:K_q\to K_{q-1}$,
$s^i_q:K_q\to K_{q+1}$, $i=0,1,\dots$.
The morphisms $\partial^i_q$ are called {\it face operators},
$s^i_q$ are called {\it degeneracy operators}.
The face- and degeneracy- operators satisfy the
following axioms:

\begin{itemize}
\item
$\partial^i_q=0$, $s^i_q=0$ if $i>q$
\item
$\partial^i_{q-1}\partial^j_q=\partial^{j-1}_{q-1}\partial^i_q$ if $i<j$
\item
$s^i_{q+1}s^j_q=s^{j+1}_{q+1}s^i_{q}$ if $i\le j$
\item
$\partial^i_{q+1}s^j_{q}=s^{j-1}_{q-1}\partial^i_{q}$ if $i<j$
\item
$\partial^i_{q+1}s^i_{q}=\partial^{i+1}_{q+1}s^i_{q}=id$ if $i\le q$
\item
$\partial^i_{q+1}s^j_{q}=s^{j}_{q-1}\partial^{i-1}_{q}$ if $i>j+1$
\end{itemize}

\medskip

\definition Given an FD-module $K$, put
$$
R(K)_q=\bigcap_{i<q}\ker(\partial^q_i:K_q\to K_{q-1}),\ \ q=0,1,\dots
$$
Clearly, $\partial^q_q(R(K)_q)\subset R(K)_{q-1}$ and
$\partial_{q-1}^{q-1}\partial^q_q|_{R(K)_q}=0$.
So $R$ is a functor from FD-modules to chain modules.
Homology groups are
$$
H_q\left(K\right)=H_q\left(R\left(K\right)\right)=
\ker\partial_q^q/\im\partial_{q+1}^{q+1}.
$$

Another way to define homology groups is as follows.
The boundary operator is
$\partial_{q}=\partial^0_{q}-\partial^1_{q}+\partial^2_{q}-\dots$,
the $q$th homology group is
$H_q=\ker\partial_q/\im\partial_{q+1}$.
Denote this functor by $R'$.

\begin{proposition} (A.~Dold, \cite{dold})
\begin{enumerate}
\item There exists a functor $R^{-1}$ from
chain modules to FD-modules such that
$RR^{-1}$ and $R^{-1}R$ are naturally equivalent
to the respective identity functors. 
\item Let $K$ be an FD-module. Chain modules
$R(K)$ and $R'(K)$ are homotopically equivalent
and therefore they have isomorphic homology groups.
\item The functor $R$ moves homotopically
equivalent FD-modules (to be defined later) to
homotopically equivalent chain modules.
\end{enumerate}
\end{proposition}

A morphism of FD-modules (or FD-morphism)
$i:K\to K'$ is a set of linear mappings
$i_q:K_q\to K'_q$ such that $i_{q-1}\partial^i=\partial^ii_q$
and $i_{q+1}s^i=s^ii_q$.
A morphism of FD-modules induces a homomorphism
of homology groups $i_{*}$.
\medskip

\example
Let $\Delta_q$ be the standard $q$-dimensional simplex:
$$
\Delta_q=\left\{(x_0,\dots,x_q)\in\R^{q+1}:x_i\ge 0,\sum x_i=1\right\}.
$$
Introduce the mappings
$$
\eps^i_q:\Delta_q\to\Delta_{q+1},\ \
\eps^i_q(x_0,\dots,x_q)=(x_0,\dots,x_{i-1},0,x_{i},\dots,x_q),
$$
$$
\eta^i_q:\Delta_q\to\Delta_{q-1},\ \
\eta^i_q(x_0,\dots,x_q)=
(x_0,\dots,x_{i-1}+x_i,x_{i+1},\dots,x_q).
$$
Suppose $X$ is a topological space,
$K_q$ is the group of its singular $q$-chains.
The elements of $K_q$ are finite sums
$\sum a_kf_k$, where $f_k:\Delta_q\to X$ are
continuous maps.
Now introduce face- and degeneracy-operators:
$$
\partial^i_q(f)=f\circ\eps^i_{q-1},
$$
$$
s^i_q(f)=f\circ\eta^i_{q+1}.
$$
We have constructed an FD-module $K(X)$ 
of a topological space $X$.

Let $x_0$ be a distinguished point of the space $X$.
All simplices $\{f:\Delta_q\to x_0\}$ form an FD-submodule
$pt(X)\subset K(X)$. The quotient $K(X)/pt(X)$ is called
the reduced FD-module of a topological space $X$,
its homology $\tilde H_*(X)$ is the reduced homology of $X$.
It is well known that
$\tilde H_q(X)\cong H_q(X)$ if $q>0$, 
$\dim H_0(X)=\dim\tilde H_0(X)+1$.

The continuous map $f:X\to Y$ induces an FD-morphism
$f_{\#}:K(X)\to K(Y)$. The map $f:(X,x_0)\to(Y,y_0)$ of pointed
spaces induces a morphism of reduced FD-modules.
\medskip

Now we define a direct product of FD-modules.

\definition
Let $K'$ and $K''$ be FD-modules. Put
$$
(K'\times K'')_q=K'_q\otimes K''_q,
$$
$$
\partial^i_{q}(a'_q\otimes
a''_q)=\partial^i_{q}(a'_q)\otimes\partial^i_{q}(a''_q),
\ s^i_{q}(a'_q\otimes a''_q)=s^i_{q}(a'_q)\otimes s^i_{q}(a''_q).
$$

Let $X'$ and $X''$ be topological spaces.
Clearly,
$K(X'\times X'')\cong K(X')\times K(X'')$. Indeed,
let $\alpha:X'\times X''\to X'$ and
$\beta:X'\times X''\to X''$ be the projections
of the Cartesian product onto each factor. Then
$\alpha_{\#}\times\beta_{\#}:K(X'\times X'')\to K(X')\times K(X'')$ 
is an isomorphism of FD-modules.
\medskip

Now define an FD-module homotopy.
First we construct an FD-module $K(n)$ that is
usually called {\it an FD-module of a standard $n$-simplex}.
The group $K(n)_q$ is free generated by all 
sets of integers $(w_0,w_1,\dots,w_q)$
such that $0\le w_0\le\dots\le w_q\le n$.
Face- and degeneracy- operators are given by the formulas
$$
\partial^i_q(w_0,\dots,w_q)=(w_0,\dots,w_{i-1},w_{i+1},\dots,w_q),
$$
$$
s^i_q(w_0,\dots,w_q)=(w_0,\dots,w_{i},w_{i},\dots,w_q).
$$
Let $e_i=(i)\in K(n)_1$ be an $i$th vertex.

\definition
Suppose $K'$ and $K''$ are FD-modules,
$F^0,F^1:K'\to K''$ are FD-morphisms.
{\it An FD-homotopy} is an FD-morphism 
$\Theta:K(1)\times K'\to K''$ such that
$$
\Theta\left((s^0_0)^qe_i\otimes a_q\right)=F^i(a_q),\
a_q\in K'_q,\ i=0,1.
$$
If there exists such homotopy $\Theta$, then 
FD-morphisms $F^0$ and $F^1$ are called homotopic,
$F^0\sim F^1$.
If there exist FD-morphisms
$F^+:K'\to K''$ and $F^-:K''\to K'$
such that
$F^+F^-\sim\id$ and $F^-F^+\sim\id$,
then FD-modules $K'$ and $K''$ are called
homotopically equivalent.
\medskip

A.~Dold has proved the following fact in \cite{dold}:

\begin{theorem}
Two FD-modules are homotopically equivalent if
and only if they have isomorphic homology groups.
\end{theorem}
\medskip

Let $T$ be a functor on the category of FD-modules.

\definition
The functor $T$ preserves homotopy if $T(F^0)\sim T(F^1)$
for each pair $F^0$, $F^1$ of homotopic FD-morphisms.

\begin{theorem} (A.~Dold, \cite{dold}) 
Suppose the functor $T$ preserves homotopy.
Then
$$
H_*(K')\cong H_*(K'')\Rightarrow H_*(T(K'))\cong H_*(T(K'')),
$$
for each pair of FD-modules $K'$, $K''$.
\end{theorem}
\medskip

Let $t$ be a functor on the category of $\Z_2$-modules.
Suppose $K$ is an FD-module. Put
$$
T(K)_q=t(K_q),\ \
F^i_{q}=t(\partial^i_{q}),\ \
D^i_{q}=t(s^i_{q}).
$$
Thus $T(K)$ is an FD-module with face-operators $F^i_q$ 
and degeneracy-operators $D^i_q$.
We say that the functor $T$ is a prolongation of the functor $t$.

The $p$th direct product $K^{\times p}$ 
of an FD-module $K$
is a simplest example of a prolongation.  

\begin{theorem} (A.~Dold, \cite{dold}) 
Suppose the functor $T$ is a prolongation
of the functor $t$. Then $T$ preserves homotopy.
\end{theorem}
\medskip

Now we can prove the following statement:

\begin{lemma}
Let $M_1$ and $M_2$ be polytopes.
Then
$H_*(M_1;\Z_2)\cong H_*(M_2;\Z_2)$ implies
$$
H_*(M_1\times M_1\times M_1/D_3,\Delta_{M_1};\Z_2)
\cong H_*(M_2\times M_2\times M_2/D_3,\Delta_{M_2};\Z_2).
$$
As above, the dihedral group $D_3$ acts on $M^{\times 3}$ by
permutations of factors,
$\Delta_M=\left\{(x,x,y)\right\}\subset M^{\times 3}/D_3$.
\end{lemma}

\beginproof
Suppose $M$ is a triangulated subset of Euclidean space,
$x_0\in M$ is a distinguished point.
Thus $(x_0,x_0,x_0)$ is a
distinguished point in the Cartesian cube $M^{\times 3}$.
Let $X=RD^3(M)$
be the reduced dihedric cube of $M$.

Let us introduce the following notations:
$$
K=K(M\times M\times M)/pt(M\times M\times M)\cong
\left(K(M)/pt(M)\right)^{\times 3},
$$
$$
K'=K(X)/pt(X).
$$

The projection $p:M^{\times 3}\to X$ induces
an FD-morphism $p_{\#}:K\to K'$.

By $K_{sym}$ denote the FD-submodule of $K$
generated by all chains of the form
$c_1\otimes c_2\otimes c_3+c_2\otimes c_3\otimes c_1$, 
$c_1\otimes c_2\otimes c_3+c_2\otimes c_1\otimes c_3$,
and $c_1\otimes c_1\otimes c_3$. 
Let $u:K_{sym}\to K$ be the inclusion,
$\varphi: K\to K/K_{sym}$ be the quotient FD-morphism.

It follows from  A.~Dold's theorem that 
the homology $H_*(K/K_{sym})$
can be computed if we know only the homology $H_*(M)$.
Let us prove that $H_*(K/K_{sym})\cong H_*(K')$, this will
complete the proof of the lemma.

First we show that
$\alpha=p_{\#}\circ\varphi^{-1}$ 
is an FD-morphism.
We need to check that
$u(K_{sym})\subset\ker p_{\#}$.

Let $c=f\otimes g\otimes h$, $c'=g\otimes h\otimes f$ 
be two elements of $K_q$.
Here $f,g,h:\Delta_q\to M$ are continuous maps of
the standard $q$-simplex to the space $M$.
Evidently, $p_{\#}(c)=p_{\#}(c')$, hence
$c+c'\in\ker p_{\#}$. 
Now suppose
$c=(\sum f_i)\otimes(\sum g_j)\otimes(\sum h_k)$,
$c'=(\sum g_j)\otimes(\sum h_k)\otimes(\sum f_i)$, 
where $f_i,g_j,h_k:\Delta_q\to M$. 
Then $c=\sum_{i,j,k}f_i\otimes g_j\otimes h_k$ and
$c'=\sum_{i,j,k}g_j\otimes h_k\otimes f_i$.
We see again that $p_{\#}(c+c')=0$.
Similarly, 
$c_1\otimes c_2\otimes c_3+c_2\otimes c_1\otimes c_3\in\ker p_{\#}$.

Suppose $c=f\otimes f\otimes g$. Then $p_{\#}(c)$ maps
the simplex $\Delta_q$ to the distinguished point,
hence $p_{\#}(c)=0\in K'$.
Now suppose $c=(\sum f_i)\otimes 
(\sum f_i)\otimes (\sum g_j)=
\sum_{i,j}f_i\otimes f_i\otimes g_j$.
We see that $p_{\#}(c)=0$.

Let us show that the FD-morphism $\alpha$ induces an 
isomorphism of homology groups.
We will construct a triangulation $\tau$
of the space $M^{\times 3}$.
The idea is as follows.
This triangulation will be coherent
with the filtration
$$
\left\{(x,x,x)\right\}=D^{0}\subset
\left\{(x,x,y),(x,y,x),(y,x,x)\right\}=D\subset
M\times M\times M
$$
and invariant under the action of the dihedric group $D_3$.
Since $\tau$ is invariant up to the action of $D_3$,
$\tau$ induces a triangulation $\tau/D_3$ of
the space $D^3(M)=M^{\times 3}$.
Since the diagonal $D$ is triangulated, 
$\tau/D_3$ induces a cell decomposition $p(\tau)$
of the space $D^3(M)/D_3$.
All affine maps $f:\Delta_q\to M^{\times 3}$ such that
$f$ moves vertices of $\Delta_q$ to vertices of
a simplex $\delta\in\tau$ generate an FD-submodule
$S\subset K$. The image $p_{\#}(S)=S'\subset K'$
is an FD-submodule of $K'$.
By $S_{sym}$ we denote the intersection 
$K_{sym}\cap S$. 
There are inclusions $i:S_{sym}\to K_{sym}$ and $j:S\to K$,
hence there is an inclusion $k:S/S_{sym}\to K/K_{sym}$.
We show all this FD-modules and morphisms on 
the following commutative diagram (its columns are exact):

$$
\Diagram
      &              & 0         &        & 0 \\
      &              & \dTo      &        & \dTo  \\
S'    &              & K_{sym}   & \lTo^i & S_{sym} \\
\dTo>l &             & \dTo>u    &        & \dTo  \\
K'    & \lTo^{p_\#}  & K         & \lTo^j & S          \\
      & \luTo^\alpha & \dTo>\varphi &     & \dTo  \\
      &              & K/K_{sym} & \lTo^k & S/S_{sym}  \\
      &              & \dTo      &        & \dTo  \\
      &              & 0         &        & 0 \\
\endDiagram
$$

In is well know that the inclusions $j$ and $l$ induce isomorphisms
of homology groups. It is easy to see that $i_*$ 
is an isomorphism too. Indeed, suppose we have a singular
cycle $\zeta\in\ker\partial_q\subset(K_{sym})_q$.
Its $0$-skeleton consists of the points on the diagonal
and pairs of symmetric points.
Thus we can construct a homotopy $F_t$ of its $0$-skeleton such that
the image of $F_t$ belongs to $K_{sym}$ for each $t$, and
the image of $F_1$ belongs to $S_{sym}$.
Then we can extend $F_t$ to $1$-skeleton, and so on.
Now two short exact sequences of FD-modules give us two
long sequences of homology groups.
By the well known $5$-lemma, it follows that $k_*$ is
an isomorphism.

Now we compute $\ker p_{\#}$.
Suppose $c\in K$, $p_{\#}(c)=0$. Clearly,
$$
c=\sum_m \left(
f_{1m}\otimes f_{2m}\otimes f_{3m}+
f_{1m}'\otimes f_{2m}'\otimes f_{3m}'
\right)
+\sum_n F_n,
$$
where $F_n:\Delta_q\to D$ are simplices in the diagonal, 
$p_{\#}(f_{1m}\otimes f_{2m}\otimes f_{3m})=
p_{\#}(f_{1m}'\otimes f_{2m}'\otimes f_{3m}')$. 
It means that
$$
f_{1m},f_{2m},f_{3m},f_{1m}',f_{2m}',f_{3m}':\Delta_q\to M
$$
and for each point $P\in\Delta_q$ there exists a permutation
$\sigma\in S_3$ such that
$\left(f_{1m}'(P),f_{2m}'(P),f_{3m}'(P)\right)=
\left(f_{\sigma(1)m}(p),f_{\sigma(2)m}(p),f_{\sigma(3)m}(p)\right)$.
Note that this permutation is the same 
for all $p\in\Delta_q$ if the image 
$(f_{1m}'\otimes f_{2m}'\otimes f_{3m}')(\Delta_q)$ does not
intersect the diagonal.
In this case 
$f_{1m}\otimes f_{2m}\otimes f_{3m}+
f_{1m}'\otimes f_{2m}'\otimes f_{3m}'\in\ker\varphi$.

Similarly, for each $F_n$ we have
$F_n=g_{1n}\otimes g_{2n}\otimes g_{3n})$, where 
$$
g_{1n},g_{2n},g_{3n}:\Delta_q\to M
$$
and for each $P\in\Delta_q$ there exists a pair of indices $i$,
$j$ such that $g_{in}(P)=g_{jn}(P)$.
If the image $F_n(\Delta_q)$ does not intersect
the small diagonal $D^0$, then this pair is the
same for all $p\in\Delta_q$.
As before, $F_n\in\ker\varphi$.

Now we see that all elements of $\ker\alpha$
have either simplices whose images intersect
the diagonal or simplices whose images lie in
the diagonal and intersect the small diagonal.
Clearly, it follows that $\ker\alpha\cap\im k=0$.
Evidently, $\im\alpha\circ k = \im l$.
We obtain that $\alpha\circ k:S/S_{sym}\to l(S')$
is an FD-isomorphism.
We see that
$l_*$, $k_*$, $l_*^{-1}\circ(\alpha\circ k)_*$
are isomorphisms.
Consequently,
$\alpha_*=l_*\circ l_*^{-1}\circ(\alpha\circ k)_* \circ k_*^{-1}$
is an isomorphism. 

Now let us construct a triangulation 
of the space $M^{\times 3}$.
Suppose the simplices $L_q^i$ form a triangulation
of the space $M$ (here $q$ is the dimension, $i$ is an index).
Then the Cartesian cube $M\times M\times M$ is divided into
sets $L_q\times L_q\times L_q$,
$L_q\times L_q\times L_r$, and $L_p\times L_q\times L_r$
(here $p$, $q$, $r$ are all distinct).

First consider $L_p\times L_q\times L_r$.
By $\tau_{pqr}$ denote
any triangulation for
$L_p\times L_q\times L_r$, $p<q<r$. 
Any triple $p,q,r$ is a permutation
of $p'<q'<r'$, so a triangulation
for $L_p\times L_q\times L_r$ is
obtained by the action of $D_3$
from $\tau_{p'q'r'}$.

Now consider $L_q\times L_q\times L_q$.
There is a filtration
$$
\left\{(x,x,x)\right\}=D^{0}\subset
\left\{(x,x,y),(x,y,x),(y,x,x)\right\}=D\subset
L_q\times L_q\times L_q.
$$
We triangulate $L_q\times L_q\times L_q$ according to this filtration
and action of the group $D_3$.
Suppose $(x,y,z)\in\Delta_q\times\Delta_q\times\Delta_q$.
Let the faces of the simplex 
$\Delta_q$ have indices $0,1,\dots,q$.
By $A^a$ denote the distance from $x$ to $a$th face,
by $B^b$ denote the distance from $y$ to $b$th face,
by $C^c$ denote the distance from $z$ to $c$th face,
$1\le a,b,c\le q$.

Suppose $0\le k\le 2q$.
Let $k_1$, $k_2$, $k_3$ be non-negative integers
such that $k_1+k_2+k_3=k$.
Suppose $\vec a=(a_1,\dots, a_{k_1+k_2})$,
$\vec b=(b_1,\dots, b_{k_1+k_3})$,
$\vec c=(c_1,\dots, c_{k_2+k_3})$, where
$1\le a_1,\dots, a_{k_1+k_2},b_1,\dots, b_{k_1+k_3},
c_1,\dots, c_{k_2+k_3}\le q$.
By $E_{\vec a,\vec b,\vec c}$ denote
the polyhedron given by equations
$$
A^{a_1}=B^{b_1},\dots,A^{a_{k_1}}=B^{b_{k_1}},
$$
$$
A^{a_{k_1+1}}=C^{c_1},\dots,A^{a_{k_1+k_2}}=C^{c_{k_2}},
$$
$$
B^{b_{k_1+1}}=C^{k_2+1},\dots,B^{b_{k_1+k_3}}=C^{c_{k_2+k_3}}.
$$
Let $Q_k$ be the union of all such polyhedra of codimension $k$:
$$
Q_k=\bigcup_{\vec a,\vec b,\vec c:\ 
\codim E_{\vec a,\vec b,\vec c}=k}
E_{\vec a,\vec b,\vec c}
$$
Then the filtration $Q_{2q}\subset Q_{2q-1}\subset\dots\subset Q_{0}$
satisfies the following conditions:

\begin{itemize}
\item $D^0=Q_{2q}$.
\item $D\subset Q_q$
\item $\Delta_q\times\Delta_q\times\Delta_q=Q_0$
\item The filtration is invariant under the action of $D_3$.
\item The set $Q_{k+1}$ divides the set $Q_k$ into pieces.
The group $D_3$ either transposes a whole piece or does not
move each point of a piece.
\end{itemize}

Now $Q_{2q}$ is triangulated. According to the last condition
the triangulation can be extended from 
the union of polyhedrons $Q_k$ to 
the union of polyhedrons $Q_{k+1}$.

Now consider $L_q\times L_q\times L_r$.
We must construct a triangulation
that is invariant up to
$D_2$-action by permuting
two $L_q$ in the Cartesian product
and coherent with the filtration
$$
\left\{(x,x,y)\right\}\subset L_q\times L_q\times L_r.
$$
First we construct a triangulation
for $L_q\times L_q$ such that it
is coherent with the filtration
$$
\left\{(x,x)\right\}\subset L_q\times L_q
$$
and invariant under the action of the group $D_2$.
We can do it in similar way as for $L_q\times L_q\times L_q$.
By $\sk_j(L_r)$ we denote a union of
all $j$-faces of the simplex $L_r$.
Now we see that for 
$L_q\times L_q\times\sk_0(L_r)$
the triangulation is constructed.
Evidently, we can extend this triangulation
from $L_q\times L_q\times\sk_j(L_r)$
to $L_q\times L_q\times\sk_{j+1}(L_r)$.

This completes the proof.
\proofend

\begin{remark}
In the same way, if
$H_*(X_1;\Z_2)\cong H_*(X_2;\Z_2)$, then
$$
H_*\left(X_1\times X_1/\Z_2,\left\{(x,x)\right\};\Z_2\right)
\cong
H_*\left(X_2\times X_2/\Z_2,\left\{(x,x)\right\};\Z_2\right).
$$
\end{remark}

\section{Proof of P.~Pushkar's theorem}

Let $M$ be an $m$-dimensional smooth closed
manifold embedded in Euclidean 
space $\R^n$, $k_q=\dim H_q(M;\Z_2)$.
We need to calculate the sum
$$
\sum\dim H_q(M^{\times 2}/\Z_2,\Delta;\Z_2).
$$
From the results of the previous section, it follows that
this sum equals
$$
\sum\dim H_q(X^{\times 2}/\Z_2,\Delta;\Z_2),
$$
where $X$ is a bouquet of spheres
$$
X=S^m\vee S^{m-1}_{1}\vee \dots\vee S^{m-1}_{k_{m-1}}
\vee\dots\vee S^{1}_{1}\vee\dots\vee S^{1}_{k_{1}}.
$$

First let us prove the following

\begin{lemma}
$$
\sum\dim H_q(S^m\times S^m/\Z_2,\Delta;\Z_2)=m+1,
$$
$$
H_*(S^m\times S^m/\Z_2,\Delta;\Z_2)=\{0,\dots,0,\Z_2,\dots,\Z_2\}.
$$
\end{lemma}

\beginproof
We construct a cell decomposition
of the reduced Cartesian square
$S^m\times S^m/\tilde\Delta$.
(Let us recall that
$\tilde\Delta=\left\{(x,x)\right\}$).
We consider $S^m$ as $m$-dimensional cube
$[0,1]^{\times m}$ with the boundary contracted to a point.
So $S^m\times S^m$ is a Cartesian product
of two such cubes. 
Let $\varphi_1,\dots,\varphi_m$ be the coordinates
on the first cube,
$\psi_1,\dots,\psi_m$ be the coordinates
on the second cube.

The cells are
$$
e^q_i=\{\varphi_1*\psi_1,\dots,
\varphi_m*\psi_m\},
$$
where $*$ is one of the signs $>$, $<$, $=$.
Of course, not all of $*$ are $=$,
since we have contracted the diagonal $\tilde\Delta$.
The codimension of a cell is $\#(=)$.
There are two $m$-dimensional cells:
$$
e^m_1=\{\varphi_1=0,\dots,\varphi_m=0\},
$$
$$
e^m_2=\{\psi_1=0,\dots,\psi_m=0\}.
$$
Evidently, they do not belong to the boundaries
of $(m+1)$-dimensional cells.

This decomposition is invariant under the $D_2$-action,
hence it induces a cell decomposition of
the reduced dihedric square $RD^2(S^m)$.\
Evidently, the chain complex for $RD^2(S^m)$
is a direct sum of two subcomplexes.
One of them is generated by all 
$e^q_i=\{\varphi_1*\psi_1,\dots,
\varphi_m*\psi_m\}$, let us denote it
by $b(RD^2(S^m))$,
and another consists of one element
$e^m=\{\varphi_1=0,\dots,\varphi_m=0\}$,
denote it by $s(RD^2(S^m))$.

Obviously,
$H_m(S^m\times S^m/\Z_2,\Delta)\cong\Z_2$: there is only one
cell of dimension $m$, its boundary is $0$.

Consider the Smith exact sequence of the Cartesian
square $S^m\times S^m$ and the group $\Z_2$ acting on it
(see \cite{viro}): 

$$
\ldots
\longrightarrow
H_q(S^m\times S^m/\Z_2,\Delta)\oplus H_q(S^m)
\longrightarrow
H_q(S^m\times S^m)
\longrightarrow
$$
$$
H_q(S^m\times S^m/\Z_2,\Delta)
\longrightarrow
H_{q-1}(S^m\times S^m/\Z_2,\Delta)\oplus H_{q-1}(S^m)
\longrightarrow
\ldots
$$

\bigskip

Since $H_q(S^m\times S^m)\cong H_q(S^m)\cong 0$ if $q\neq 0,m,2m$, 
we have
$$
0
\longrightarrow
H_q(S^m\times S^m/\Z_2,\Delta)
\longrightarrow
H_{q-1}(S^m\times S^m/\Z_2,\Delta)
\longrightarrow
0,
$$ 
where $q=m+2,\dots, 2m-1$.

Thus we see that 
$H_{m+1}(S^m\times S^m/\Z_2,\Delta)\cong\dots\cong
H_{2m-1}(S^m\times S^m/\Z_2,\Delta)$.

Consider the end of the Smith sequence:
$$
0
\longrightarrow
H_{m+1}(S^m\times S^m/\Z_2,\Delta)
\longrightarrow
H_{m}(S^m\times S^m/\Z_2,\Delta)\oplus\Z_2
\longrightarrow
$$
$$
\Z_2\oplus\Z_2
\longrightarrow
H_{m}(S^m\times S^m/\Z_2,\Delta)
\longrightarrow
0.
$$

Since $H_{m}(S^m\times S^m/\Z_2,\Delta)\cong\Z_2$, we see
that $H_{m+1}(S^m\times S^m/\Z_2,\Delta)\cong\Z_2$.
Thus we have
$$
H_{m}(S^m\times S^m/\Z_2,\Delta)\cong\dots\cong
H_{2m-1}(S^m\times S^m/\Z_2,\Delta)\cong\Z_2.
$$

Now consider the beginning of the Smith sequence:
$$
0
\longrightarrow
H_{2m}(S^m\times S^m/\Z_2,\Delta)
\longrightarrow
\Z_2
\longrightarrow
H_{2m}(S^m\times S^m/\Z_2,\Delta)
\longrightarrow
$$
$$
H_{2m-1}(S^m\times S^m/\Z_2,\Delta)
\longrightarrow
0.
$$
We already know that $H_{2m-1}(S^m\times S^m/\Z_2,\Delta)\cong\Z_2$.
Computing the Euler characteristic we obtain that
$H_{2m}(S^m\times S^m/\Z_2,\Delta)\cong\Z_2$.
This completes the proof.
\proofend

Now let us prove P.~Pushkar's estimate.

\smallskip
\begin{lemma}
Let $M$ be an $m$-dimensional manifold, 
$k_i=\dim H_i(M;\Z_2)$,
$i=0,1,\dots,m$, $B=\sum k_i$. Then
$$
\sum_{i=1}^mik_i=\frac{mB}{2}.
$$
\end{lemma}
\beginproof
Using Poincar\'e duality, we have
$$
\sum_{i=1}^m ik_i=
\sum_{i=0}^m ik_i=
\frac{1}{2}\sum_{i=0}^m\left(ik_i+\left(m-i\right)k_{m-i}\right)=
\frac{1}{2}\sum_{i=0}^m\left(ik_i+\left(m-i\right)k_{i}\right)=
\frac{mB}{2}.
$$
\proofend

\begin{theorem}
Let $B=\sum\dim H_q(M;\Z_2)$. Then
$$
BT_2(M)\ge\frac{B^2+(m-1)B}{2}.
$$
\end{theorem}

\beginproof

Let us recall that $X$ 
is a bouquet of spheres such that
$H_*(X;\Z_2)\cong H_*(M;\Z_2)$.
Taking into account all arguments above, we see that
it is enough to prove that
$$
(*)\ \ \ 
\sum\dim H_q(X^{\times 2}/\Z_2,\Delta;\Z_2)=\frac{B^2+(m-1)B}{2}.
$$

We construct a cell decomposition of the space
$RD^2(X)$. This space is the union of its subsets:
$$
RD^2(X)=\bigcup_{{1\le p\le m}\atop{1\le i\le k_p}} A^p_i \cup
\bigcup_{{1\le p\le m}\atop{1\le i<j\le k_p}} B^p_{ij} \cup
\bigcup_{{1\le p<q\le m}\atop{1\le i\le k_p,1\le j\le k_q}} C^{pq}_{ij},
$$
where

\begin{itemize}
\item $A^p_i$ is the set of points
$(x,x')\in RD^2(X)$ such that $x,x'\in S^p_i$.
Thus $A^p_i$ is homeomorphic to $RD^2(S^p)$.

\item $B^p_{ij}$ is the set of points
$(x_1,x_2)\in RD^2(X)$ such that $x_1\in S^p_i$,
$x_2\in S^p_j$.
Thus $B^p_{ij}$ is homeomorphic to $S^p\times S^p$.

\item $C^{pq}_{ij}$ is the set of points
$(x,y)\in RD^2(X)$ such that $x\in S^p_i$,
$y\in S^q_j$.
Thus $C^{pq}_{ij}$ is homeomorphic to $S^p\times S^q$.
\end{itemize}

For each $A^p_i$ we have constructed a cell
decomposition earlier.
For each $B^p_{ij}$, $C^{pq}_{ij}$
we have a standard decomposition:
$$
S^p=e^0\cup e^p,\ \
S^q=e^0\cup e^q,\ \
S^p\times S^q=e^0\times e^0\cup
e^0\times e^q\cup
e^p\times e^0\cup
e^p\times e^q
$$
with all boundaries equal to $0$.
By $s(B^p_{ij})$ 
we denote a chain complex generated
by one cell $e^p\times e^p$,
by $s(C^{pq}_{ij})$
we denote a chain complex generated
by one cell $e^p\times e^q$.

Notice that these decompositions
coincide on intersections of the sets $A^p_i$,
$B^p_{ij}$, $C^{pq}_{ij}$.

Now we see that a chain complex for $RD^2(X)$
is a direct sum of its subcomplexes:
$$
\bigoplus_{{1\le p\le m}\atop{1\le i\le k_p}}
\left(b(A^p_i)\oplus s(A^p_i)\right)
\bigoplus_{{1\le p\le m}\atop{1\le i<j\le k_p}} s(B^p_{ij})
\bigoplus_{{1\le p<q\le m}\atop{1\le i\le k_p,1\le j\le k_q}} s(C^{pq}_{ij})
$$

Now we compute contributions to the sum of Betti
numbers $(*)$ for each of these direct summands.

\begin{itemize}

\item{$b(A^p_i)\oplus s(A^p_i)$: The sum of Betti numbers of
each $A^p_i$ is $p+1$. So we have a contribution
$\sum_{i=1}^m (i+1)k_i$.

Let us transform this expression using the previous lemma.
$$
\sum_{p=1}^m (p+1)k_p=
B-1+\frac{mB}{2}.
$$
}

\item{$s(B^p_{ij})$: 
Each of $s(B^p_{ij})$ gives us $1$ to the sum
of Betti numbers $(*)$, so totally we have
$$
\frac{1}{2}\sum_{p=1}^m k_p(k_p-1).
$$
}

\item{$s(C^{pq}_{ij})$:
$$
\sum_{p<q}^m k_pk_q.
$$
}
\end{itemize}

Thus we have
$$
\sum\dim H_q(X^{\times 2}/\Z_2,\Delta;\Z_2)=
B-1+\frac{mB}{2}+
\frac{1}{2}\sum_{p=1}^m k_p(k_p-1)+
\sum_{p<q}^m k_pk_q=
$$
$$
B-1+\frac{mB}{2}+
\frac{1}{2}\left(\sum_{p=1}^m k_p\right)^2-
\frac{1}{2}\sum_{p=1}^m k_p=
\frac{(B-1)^2+(m+1)B-1}{2}=
\frac{B^2+(m-1)B}{2}.
$$
This completes the proof.
\proofend

\section{Estimate for the number of $3$-periodical trajectories}

As above, a smooth closed manifold $M$ is embedded
in Euclidean space $\R^n$. 
Let $k_q=\dim H_q(M;\Z_2)$.
We need to compute the sum
$$
\sum\dim H_q(M^{\times 3}/S_3,\Delta;\Z_2).
$$
It equals to
$$
\sum\dim\tilde H_q\left(RD^3(X);\Z_2\right),
$$
where $X$ is a bouquet of spheres:
$$
X=S^m\vee S^{m-1}_{1}\vee \dots\vee S^{m-1}_{k_{m-1}}
\vee\dots\vee S^{1}_{1}\vee\dots\vee S^{1}_{k_{1}},
$$
We know that (see \cite{farber},\cite{farber2})
$$
\sum\dim\tilde H_q\left(RD^3(S^m);\Z_2\right) =2m.
$$

First let us prove the following
\begin{lemma}
Let $M$ be $m$-dimensional manifold,
 $k_i=\dim H_i(M;\Z_2)$,
$i=0,1,\dots,m$, $B=\sum k_i$. Then
$$
\sum_{1\le i\le m}ik_i^2+\sum_{1\le i\neq j\le m}ik_ik_j=
\frac{mB^2-mB}{2}.
$$
\end{lemma}
\beginproof
Taking into account
Poincar\'e duality $k_i=k_{m-i}$, we have
$$
\sum_{1\le i\le m}ik_i^2+\sum_{1\le i\neq j\le m}ik_ik_j=
$$
$$
=\sum_{0\le i\le m}ik_i^2+\sum_{0\le i\neq j\le m}ik_ik_j-
\sum_{1\le i\le m, j=0}ik_ik_j=
$$
$$
=\sum_{0\le i,j\le m}ik_ik_j-
\sum_{0\le i\le m}ik_i=
$$
$$
=\frac{1}{2}\sum_{0\le i,j\le m}(i+m-i)k_ik_j-
\frac{1}{2}\sum_{0\le i\le m}(i+m-i)k_i=
$$
$$
=\frac{mB^2-mB}{2}.
$$
\proofend

\medskip
\begin{theorem}
Suppose $B=\sum_{i=0}^m k_i$. Then we have an estimate
for the number of $3$-periodical trajectories in $M$
$$
BT_3(M)\ge\frac{B^3+3(m-1)B^2+2B}{6}.
$$
\end{theorem}

\beginproof

Let us construct a cell decomposition of
the reduced dihedric cube of $X$.
The space $RD^3(X)$ is a union of its
subsets:
$$
RD^3(X)=\bigcup_{{1\le p\le m}\atop{1\le i\le k_p}} A^p_i \cup
\bigcup_{{1\le p\le m}\atop{1\le i\neq j\le k_p}} B^p_{ij} \cup
\bigcup_{{1\le p\le m}\atop{1\le i<j<k\le k_p}} C^{p}_{ijk} \cup
$$
$$
\bigcup_{{1\le p\neq q\le m}
\atop{1\le i\le k_p,1\le j\le k_q}} D^{pq}_{ij} \cup
\bigcup_{{1\le p\neq q\le m}
\atop{1\le i<j\le k_p,1\le k\le k_q}} E^{pq}_{ijk} \cup
\bigcup_{{1\le p<q<r\le m}
\atop{1\le i\le k_p,1\le j\le k_q,1\le k\le k_r}} F^{pqr}_{ijk}
$$
Here
\begin{itemize}
\item $A^p_i$ is the set of points
$(x,x',x'')\in RD^3(X)$ such that $x,x',x''\in S^p_i$.
Thus it is homeomorphic to 
$$
RD^3(S^p).
$$

\item $B^p_{ij}$ is the set of points
$(x_1,x_1',x_2)\in RD^3(X)$ such that $x_1,x_1'\in S^p_i$,
$x_2\in S^p_j$.
Thus it is homeomorphic to
$$
RD^2(S^p)\times S^p/\{(x,O,O)\},
$$ 
where $O$
is the distinguished point of the bouquet $X$.

\item $C^{p}_{ijk}$ is the set of points
$(x_1,x_2,x_3)\in RD^3(X)$ such that $x_1\in S^p_i$,
$x_2\in S^p_j$, $x_3\in S^p_k$.
Thus it is homeomorphic to
$$
S^p\times S^p\times S^p/\tilde\Delta_{S^p}.
$$

\item $D^{pq}_{ij}$ is the set of points
$(x,x',y)\in RD^3(X)$ such that $x,x'\in S^p_i$,
$y\in S^q_j$.
Thus it is homeomorphic to
$$
RD^2(S^p)\times S^q/\{(x,O,O)\}.
$$

\item $E^{pq}_{ijk}$ is the set of points
$(x_1,x_2,y)\in RD^3(X)$ such that $x_1\in S^p_i$,
$x_2\in S^p_j$, $y\in S^q_k$.
Thus it is homeomorphic to
$$
S^p\times S^p\times S^q/
\left\{(x_1,O,O) \right\},
\left\{(O,x_2,O) \right\},
\left\{(O,O,y) \right\}.
$$

\item $F^{pqr}_{ijk}$ is the set of points
$(x,y,z)\in RD^3(X)$ such that $x\in S^p_i$,
$y\in S^q_j$, $z\in S^r_k$.
Thus it is homeomorphic to
$$
S^p\times S^q\times S^r/
\left\{(x,O,O) \right\},
\left\{(O,y,O) \right\},
\left\{(O,O,z) \right\}.
$$
\end{itemize}

Let us construct a cell decomposition for
each $A^p_i=RD^3(S^p)$ like we have done
it for $RD^2(S^p)$.
We consider a $p$-sphere as
a $p$-dimensional cube $[0,1]^{\times p}$ with
the boundary contracted to a point.
Let $\varphi^1_i$, 
$\varphi^2_j$, $\varphi^3_k$
be the coordinates on three such cubes.
There are the following cells.
First
$$
U^d_{\beta}=\{\varphi_1^{i_1}*\varphi_1^{j_1}*\varphi_1^{k_1},\dots,
\varphi_m^{i_m}*\varphi_m^{j_m}*\varphi_m^{k_m}\},
$$
where $*$ is one of the signs $<$ or $=$,
$(i_{\alpha},j_{\alpha},k_{\alpha})$ 
is a permutation of $1$, $2$, $3$. 
Secondly
$$
V^d_{\beta}=\{\varphi_1^1*\varphi_1^2,\dots,
\varphi_m^1*\varphi_m^2,\vec\varphi^3=0\},
$$
where $*$ is one of the signs $<$, $>$, $=$.
Clearly, the cells of the kind $U$ do not
belong to the boundaries of the cells of the kind $V$.
Similarly, the cells of the kind $V$ do not
belong to the boundaries of the cells of the kind $U$.
It is evident for $p>1$. For $p=1$ all boundary operators
here vanish.

We see that a chain complex for $A^p_i$ is a direct sum
of two subcomplexes $u(A^p_i)$ and $v(A^p_i)$.

Now we construct a cell decomposition for $B^p_{ij}$
and $D^{pq}_{ij}$.
Here we have the same cells ``of the kind $V$'':
$$
V^d_{\beta}=\{\varphi_1^1*\varphi_1^2,\dots,
\varphi_m^1*\varphi_m^2,\vec\psi=0\},
$$
where $\vec\varphi^1$ and $\vec\varphi^2$
are the coordinates on $i$th $p$-sphere,
$\vec\psi$ are the coordinates on
$j$th $p$- or $q$-sphere.
Besides, there are the cells
$$
W^d_{\beta}=\{\varphi_1^1*\varphi_1^2,\dots,
\varphi_m^1*\varphi_m^2,0<\psi_k<1\},
$$
$$
\tilde W^d_{\beta}=\{\vec\varphi^1=0,0<\varphi_k^2<1,0<\psi_l<1\}.
$$

As above, the cells of each of these three kinds form
a chain subcomplex. Thus a chain complex for 
$B^p_{ij}$ is a direct sum
$v(B^p_{ij})\oplus w(B^p_{ij})\oplus\tilde w(B^p_{ij})$
and so is one for $D^{pq}_{ij}$.

Notice that 
$$
v(B^p_{ij})=v(D^{pq}_{ij})=v(A^p_i),
$$
$$
\tilde w(B^p_{ij})=\tilde w(B^p_{ji}),\ \ 
\tilde w(D^{pq}_{ij})=\tilde w(D^{qp}_{ji}).
$$
Clearly, each $\tilde w(\dots)$ consists
of only one cell with boundary equal to $0$.

Finally, $C^p_{ijk}$, $E^{pq}_{ijk}$, and $F^{pqr}_{ijk}$
are Cartesian products of spheres, so here we have a
standard cell decomposition:
$$
S^p=e^0\cup e^p,\ \
S^q=e^0\cup e^q,\ \
S^r=e^0\cup e^r,
$$
$$
S^p\times S^q\times S^r=O\cup e^0\times e^q\times e^r \cup
e^p\times e^0\times e^r \cup
e^p\times e^q\times e^0 \cup
e^p\times e^q\times e^r.
$$
We see that $e^0\times e^q\times e^r$,
$e^p\times e^0\times e^r$, and
$e^p\times e^q\times e^0$ are cells
``of the kind $\tilde W$''.
The cell $e^p\times e^q\times e^r$ with
vanishing boundary forms a chain subcomplex.
We denote this subcomplex by
$s(C^p_{ijk})$, $s(E^{pq}_{ijk})$, or $s(F^{pqr}_{ijk})$.

Now we see that the chain complex for
$RD^3(X)$ is a direct sum of subcomplexes:
$$
\bigoplus_{{1\le p\le m}\atop{1\le i\le k_p}}
\left(
 u(A^p_i)\oplus v(A^p_i)
\right)
\bigoplus_{{1\le p\le m}\atop{1\le i\neq j\le k_p}} 
w(B^p_{ij}) 
\bigoplus_{{1\le p\le m}\atop{1\le i<j\le k_p}} 
\tilde w(B^p_{ij}) 
\bigoplus_{{1\le p\le m}\atop{1\le i<j<k\le k_p}} 
s(C^{p}_{ijk})
$$
$$
\bigoplus_{{1\le p\neq q\le m}
\atop{1\le i\le k_p,1\le j\le k_q}} 
w(D^{pq}_{ij})
\bigoplus_{{1\le p<q\le m}
\atop{1\le i\le k_p,1\le j\le k_q}} 
\tilde w(D^{pq}_{ij})
$$
$$
\bigoplus_{{1\le p\neq q\le m}
\atop{1\le i<j\le k_p,1\le k\le k_q}} s(E^{pq}_{ijk})
\bigoplus_{{1\le p<q<r\le m}
\atop{1\le i\le k_p,1\le j\le k_q,1\le k\le k_r}} s(F^{pqr}_{ijk})
$$

Now let us consider all these subcomplexes and compute
their contributions to the sum of Betti numbers
$\sum\dim\tilde H_q(RD^3(X);\Z_2)$.

\begin{itemize}

\item{$u(A^p_i)\oplus v(A^p_i)$:
$2\sum_{p=1}^m pk_p=mB$.

Here we have
$$
mB.
$$
}

\item $w(B^p_{ij})$: 
$\sum_{p=1}^m pk_p(k_p-1)$.
Indeed, for each $p$ we have $k_p(k_p-1)$ sets
$B^p_{ij}$; the cells of the kind $W$ form
the reduced dihedric square of the $p$-sphere
without the unique $p$-dimensional cell 
that forms a subcomplex $\tilde w(B^p_{ij})$.

\item $\tilde w (B^p_{ij})$: 
$\frac{1}{2}\sum_{p=1}^m k_p(k_p-1)$.
Here $\frac{1}{2}$ appears, since
$\tilde w(B^p_{ij})=\tilde w(B^p_{ji})$.

\item$w(D^{pq}_{ij})$: 
$\sum_{p\neq q} pk_pk_q$.

\item$\tilde w(D^{pq}_{ij})$: 
$\sum_{1\le p<q\le m}^m k_pk_q$.

\item{$s(C^{p}_{ijk})$: Here we have
subcomplexes consisting of one cell.
Thus we have a contribution
$$
\frac{1}{6}\sum_{p=1}^m k_p(k_p-1)(k_p-2).
$$
}

\item{$s(E^{pq}_{ijk})$: Here we have
$$
\frac{1}{2}\sum_{p\neq q}^m k_p(k_p-1)k_q.
$$
}

\item{$s(F^{pqr}_{ijk})$: Here we have
$$
\sum_{p<q<r}^m k_pk_qk_r.
$$
}
\end{itemize}

Taking into account two computational lemmas, 
we have
$$
mB+
\sum_{p=1}^m pk_p(k_p-1)+\frac{1}{2}\sum_{p=1}^m k_p(k_p-1)+
\sum_{p\neq q} pk_pk_q+\sum_{1\le p<q\le m}^m k_pk_q+
$$
$$
+\frac{1}{6}\sum_{p=1}^m k_p(k_p-1)(k_p-2)+
\frac{1}{2}\sum_{p\neq q}^m k_p(k_p-1)k_q+
\sum_{p<q<r}^m k_pk_qk_r=
$$
$$
=mB+\frac{mB^2}{2}-\frac{mB}{2}-\frac{mB}{2}+
\frac{(B-1)^2}{2}-\frac{B-1}{2}+
\frac{(B-1)^3}{6}-
\frac{(B-1)^2}{2}+\frac{B-1}{3}=
$$
$$
=\frac{1}{6}(B^3-3B^2+3B-1+3mB^2-B+1)=
$$
$$
=\frac{B^3+3(m-1)B^2-2B}{6}.
$$
This completes the proof.
\proofend


\newpage

\begin{thebibliography}{99}
\bibitem{birkhoff}
G.~Birkhoff.
Dynamical systems.
New York, 1927.
\bibitem{pushkar}
P.~Pushkar.
Periodic trajectories of billiard dynamical systems.
Ph.D. thesis. Moscow, 1998.
\bibitem{farber}
M.~Farber, S.~Tabachnikov.
Topology of cyclic configuration spaces and
periodic trajectories of multi-dimensional billiards.
arXiv:math.DG/9911226.
\bibitem{farber2}
M.~Farber, S.~Tabachnikov.
Periodic trajectories in 3-dimensional convex billiards.
arXiv:math.DG/0106049.
\bibitem{duzhin}
F.~Duzhin. On the lower bounds for the number of
periodic billiard trajectories in manifolds
embedded in Euclidean space.
arXiv:math.AT/0112178.
\bibitem{dold}
A.~Dold. Homology of symmetric products and other functors of complexes.
Annals of Mathematics, vol.~68, No.~1, July, 1958.
\bibitem{eilenberg}
S.~Eilenberg, S.~MacLane. On the groups $H(\Pi,n)$.
Annals of Mathematics, vol.~58, No.~1, July, 1953.
\bibitem{viro} 
O.~Viro, D.~Fuchs. Topology-2. In: Encyclopedia of Math. Sciences,
vol. 24.
\bibitem{griffiths}
P.~Griffiths, J.~Harris.
Principles of Algebraic Geometry.
John Wiley \& Sons, 1978.
\bibitem{fomenko}
A.~Fomenko, D.~Fuchs.
A Course in Homotopic Topology.
\bibitem{dnf}
B.~Dubrovin, S.~Novikov, A.~Fomenko.
Modern Geometry.
\end{thebibliography}
\end{document}